# General branching processes in discrete time as random trees

PETER JAGERS* and SERIK SAGITOV**

*School of Mathematical Sciences, Chalmers University of Technology, SE-412 96 Gothenburg, Sweden.* *E-mail: jagers@chalmers.se; ** serik@chalmers.se

The simple Galton–Watson process describes populations where individuals live one season and are then replaced by a random number of children. It can also be viewed as a way of generating random trees, each vertex being an individual of the family tree. This viewpoint has led to new insights and a revival of classical theory.

We show how a similar reinterpretation can shed new light on the more interesting forms of branching processes that allow repeated bearings and, thus, overlapping generations. In particular, we use the stable pedigree law to give a transparent description of a size-biased version of general branching processes in discrete time. This allows us to analyze the $x \log x$ condition for exponential growth of supercritical general processes as well as relation between simple Galton–Watson and more general branching processes.

## 1. Introduction

Branching processes are studied in at least four different traditions, those of classical branching processes, of random trees, of general (non-Markov) branching and of superprocesses, including branching diffusions and ditto random walks. The classical approach, as established in text books, is that of conventional, largely analytical, applied probability. It treats Galton–Watson (GW), birth-and-death and Markov branching processes. During the last decade-and-a-half, there has been a revival of GW processes connected with random graphs and trees, and partly inspired by computer science. All along, a more biologically or demographically motivated line of thinking has been investigating models where lifespans are not exponentially distributed or degenerate and births are admitted at several bearing ages, so as to better mimic real life. This, however, leads to a comparatively inaccessible framework and, as a consequence, to GW-style processes and trees having been independently developed. The purpose of this paper is to lay bare results from general branching processes by specializing to the case of discrete time structure, to interpret them in terms of random trees and to compare with the simple Galton–Watson case.

A GW branching process generates a random tree via a simple growth algorithm governed by a discrete probability distribution $\{p_k\}_{k=0}^{\infty}$ with $p_1 < 1$ and mean $m = \sum_{k=1}^{\infty} kp_k$. In classical parlance, this is the *reproduction* or *offspring* distribution. The algorithm starts by placing a root vertex at zero level $n = 0$. From the root, it draws $k$ branches with probability $p_k$. The ends of the $k$ branches become the tree vertices at level $n = 1$. If $k \geq 1$, then each of the $k$ vertices produces new branches, resulting in a random number of vertices at the level $n = 2$, which, in turn, give rise to branches reaching the level $n = 3$, and so on. The key assumption is that the numbers of branches stemming from different vertices are independent and identically distributed. The height of the random tree is of course the generation or discrete time (year or season) in the classical interpretation of a population with non-overlapping generations and one ancestor. If the final height of a GW tree is finite, the population ultimately dies out. The height is infinite with a positive probability if and only if reproduction is supercritical, $m > 1$.

A modified way of growing random trees, using so-called *size-biased GW* (abbreviated as $\widehat{\mathrm{GW}}$) processes, was suggested in [8]. It describes particle populations which never die out because, at any time, they contain exactly one particle with a different, *size-biased* reproduction law (for some discussion of its background, see Section 3)

$$\hat{p}_k = m^{-1} k p_k, \qquad k \geq 1. \tag{1}$$

Since its progeny never dies out, we call this particle 'immortal', for ease of reference. A child, chosen at random among the $k$ offspring of the immortal, becomes the next generation's immortal individual and the remaining $k - 1$ children are mortals whose reproduction is described by the original offspring distribution, $\{p_k\}_{k=0}^{\infty}$, the same being true for their children, and so on. $\widehat{\mathrm{GW}}$ trees are infinite and contain a distinguished path, the *immortal lineage*. The construction yields a spinal decomposition of GW processes which was used to provide alternative and illuminating proofs of the classical limit theorems for GW processes in [8] and [2]. The approach in [8] hinges on the fact that $\widehat{\mathrm{GW}}$ measures on the space of trees have a nice Radon–Nikodym derivative with respect to GW measures. (See [1] for a recent development of this approach to GW processes with a general-type space.)

In this paper, we consider branching processes with overlapping generations, called general branching processes in [3] and sometimes referred to as CMJ (Crump–Mode–Jagers) processes. For the sake of transparency, we restrict ourselves to discrete time and provide an explicit (and minimal) construction of the process. We try to describe size-biased CMJ processes simply and lucidly by using the stable pedigree law obtained in [6] and [10]. Our construction clarifies the more general case considered in [13]. The change-of-measure technique used allows us to specify the $x \log x$ condition for supercritical CMJ processes with finite mean age at childbearing in a precise way.

## 2. Multiple birth processes in discrete time

CMJ processes depict populations of individuals who can produce offspring at different ages. The key assumption is that individuals reproduce independently and according to a common reproduction law. We thus restrict ourselves to *single-type* populations in *discrete time*. The individual reproduction law is then given by a set of probabilities,

$$p_0, p_k(n_1, \ldots, n_k), \qquad 1 \le n_1 \le \cdots \le n_k < \infty, k \ge 1, \tag{2}$$

where all $n_k \in \mathbb{N} = \{1, 2, 3, \ldots\}$, $p_0$ stands for the probability of no offspring, $p_1(n_1)$ is the probability of having exactly one child who is born when the mother is $n_1$ years old, $p_2(n_1, n_2)$ is the probability that the potential mother gives birth to two daughters, first at age $n_1$ and then at age $n_2$, and so on. As in the GW case, we denote the marginal distribution of the offspring number by

$$p_k = \sum_{1 \le n_1 \le \cdots \le n_k < \infty} p_k(n_1, \ldots, n_k), \qquad k \ge 1,$$

with the mean $m = \sum_{k=1}^{\infty} k p_k$.

Reproduction yields a minimal description of individual life. It can be enriched, for example, by information about lifespan, body mass or DNA content at different ages, or other aspects of interest. The minimal setting is, however, sufficient for the basic population dynamics and we proceed to an explicit construction of the corresponding probability space for CMJ processes stemming from a single progenitor born at time zero. All possible reproduction scenarios for an individual are covered by the countable set $\Omega_0 = \bigcup_{k=0}^{\infty} \mathbb{N}_\uparrow^k$, where $\mathbb{N}_\uparrow^0 = \{0\}$ is the outcome of no offspring and $\mathbb{N}_\uparrow^k$ is the set of vectors $(n_1, \ldots, n_k)$ with $n_1 \le \cdots \le n_k$ and $n_i \in \mathbb{N}$. A discrete probability space $(\Omega_0, \mathcal{A}_0, P_0)$ for individual life is straightforwardly determined by the set of probabilities (2) via

$$P_0(\omega_0) = p_k(n_1, \ldots, n_k), \qquad \omega_0 = (n_1, \ldots, n_k).$$

An element $\omega_0 = (n_1, \ldots, n_k)$ of $\Omega_0$ determines the value of random variables like the offspring number $\nu$, $\nu(\omega_0) = k$, and the consecutive ages at bearing $\tau_i$,

$$\tau_i(\omega_0) = n_i, \qquad i = 1, \ldots, k; \qquad \tau_i = \infty, \qquad i > k. \tag{3}$$

This allows us to write

$$P_0(\nu = k, \tau_1 = n_1, \ldots, \tau_k = n_k) = p_k(n_1, \ldots, n_k)$$

for the joint distribution of bearing ages and offspring number. Putting $\tau_0 = 0$, we can refer to $\tau_\nu$ as the last time of birth. Notice that $\nu = 0$ is equivalent to $\tau_\nu = 0$.

A CMJ process realization can also be viewed as a planar tree growing from the zero level upward. In contrast to the GW picture, the height or time parameter then counts years rather than generations. A branch $\omega_0 \in \Omega_0$ of the tree represents a random individual life: if $\omega_0 = (n_1, \ldots, n_k)$, then $k$ is the branch length and $n_1, \ldots, n_k$ are the heights at which the daughter vertices appear along the branch.

The set $\mathbb{I} = \bigcup_{j=0}^{\infty} \mathbb{N}^j$ consists of all conceivable vertices in a rooted tree. Here, $\mathbb{N}^0 = \{0\}$ contains the label of the ancestor, $\mathbb{N}^1 = \{1, 2, 3, \ldots\}$ lists possible daughters in the order of appearance (ties are resolved arbitrarily), $\mathbb{N}^2$ is the set of labels for the granddaughters, and so on, with $\mathbb{N}^j = \mathbb{N}^{j-1} \times \mathbb{N}$. The usual convention is to denote the $i$th daughter of individual $x = (x_1, \ldots, x_j)$ by $xi = (x_1, \ldots, x_j, i)$.

To each potential vertex $x \in \mathbb{I}$, assign a random life history $\omega_x \in \Omega_0$, all of these independent and distributed like $\omega_0$. In strict terminology, for each $x$, let $(\Omega_x, \mathcal{A}_x)$ be a copy of $(\Omega_0, \mathcal{A}_0)$ and consider the product probability space

$$(\Omega, \mathcal{A}, P) = \prod_{x \in \mathbb{I}} (\Omega_x, \mathcal{A}_x, P_0).$$

A sequence $\omega = \{\omega_x\}_{x \in \mathbb{I}} \in \Omega$ then contains information about all possible (and maybe some impossible) branches. Only some of them are actually building the random tree, or actually get born, to speak in terms of population dynamics. Different $\omega$ may thus result in the same tree. The actual population tree can be formally defined as the set of actual branches

$$T(\omega) = \{\omega_x, x : \sigma_x < \infty\},$$

where

$$\sigma_x(\omega) = \tau_{x_1}(\omega_0) + \tau_{x_2}(\omega_{x_1}) + \cdots + \tau_{x_{j-1}}(\omega_{x_1 \cdots x_{j-2}}) + \tau_{x_j}(\omega_{x_1 \cdots x_{j-1}})$$

is the height of the vertex $x = (x_1, \ldots, x_j)$ and, demographically speaking, that individual's birthtime. In view of (3), the equality $\sigma_x = \infty$ identifies the vertices $x \in \mathbb{I}$ that do not belong to the corresponding tree. The distribution of $T$ under $P$ induces a probability measure $P(T = t)$ on the tree space (see [12] for a detailed description). This will be compared to another probability measure on the same tree space, the *size-biased* measure $\hat{P}$.

## 3. The immortals

As explained in Section 1, size-biased GW processes involve immortal individuals with a modified offspring distribution (1). The latter has a long history and describes the number of children of a mother chosen by the following indirect sampling procedure.

Consider a population of $n$ individuals who produce children independently of each other, $k$ of them with probability $p_k$. Sample one child at random and call her 'Ego' [7]. Two questions can be posed about her:

1. What is the size $S_n$ of her sibship?
2. What number is Ego in her sibship?

The second question may seem pointless since numbering is arbitrary in the present GW-type situation, but in more elaborate models like the CMJ process, this is the question of birth rank [4]. In the simple context of Galton–Watson branching, Ego's label $R_n$ will be uniformly distributed over $1, 2, \ldots, S_n$. As $n \to \infty$,

$$P(S_n = k) \to m^{-1} k p_k,$$

by the law of large numbers, and the joint distribution of sibship and rank converges to

$$\hat{p}_{k,j} = 1_{\{1 \le j \le k\}} m^{-1} p_k. \tag{4}$$

The marginal limiting distribution of rank, as $n \to \infty$, reduces to

$$\tilde{p}_j = \sum_{k=j}^{\infty} \hat{p}_{k,j} = m^{-1} \sum_{k=j}^{\infty} p_k = m^{-1} P_0(\nu \ge j). \tag{5}$$

The $\widehat{\mathrm{GW}}$ process consists of an ancestor initiating an infinite line of immortals, each reproducing and choosing the next immortal according to (5), all other individuals having offspring according to the original reproduction distribution $\{p_k\}$.

The situation for models with overlapping generations is considerably more involved. If we start from $n$ individuals, as above, and wait until all of those have completed reproduction, then there are many more individuals born than those who stem directly from the ancestors. Moreover, in a growing population, there will be more children born by young mothers than by old ones. In that supercritical case, the solution is to sample from an old general branching process, where the age distribution stabilizes, if the population avoids extinction. This was performed in [6] and [10], and the size-biased life distribution generalizing (4) was found to be given by the formula

$$\hat{p}_{k,j}(n_1, \ldots, n_k) = 1_{\{1 \le j \le k\}} \mathrm{e}^{-\alpha n_j} p_k(n_1, \ldots, n_k). \tag{6}$$

Here, $\alpha$ is the *Malthusian parameter* or *net reproduction rate*, defined by the equation

$$\sum_{n=1}^{\infty} \mathrm{e}^{-\alpha n} m_n = 1, \tag{7}$$

where $m_n$ is the mean litter number of a mortal individual at age $n \in \mathbb{N}$,

$$m_n = \mathrm{E}_0 \left( \sum_{j=1}^{\infty} 1_{\{\tau_j = n\}} \right)$$

$$= \sum_{j=1}^{\infty} P_0(\tau_j = n)$$

$$= \sum_{k=1}^{\infty} \sum_{j=1}^{k} \sum_{1 \le n_1 \le \cdots \le n_k < \infty} 1_{\{n_j = n\}} p_k(n_1, \ldots, n_k).$$

Since the sum of the litter sizes gives the total offspring number we have

$$\sum_{n=1}^{\infty} m_n = m.$$

Equation (7) has a unique solution unless $m < 1$ and $m_n$ tends to zero as $n \to \infty$ so slowly that $\sum_{n=1}^{\infty} e^{\epsilon n} m_n = \infty$, regardless of how small $\epsilon > 0$ is (as is the case for $m_n \asymp n^{-c}$, $1 < c < \infty$).

The summands in (7),

$$e^{-\alpha n} m_n = \sum_{k=1}^{\infty} \sum_{j=1}^{k} \sum_{1 \le n_1 \le \cdots \le n_k < \infty} 1_{\{n_j = n\}} e^{-\alpha n} p_k(n_1, \ldots, n_k)$$

$$= \sum_{k=1}^{\infty} \sum_{j=1}^{k} \sum_{1 \le n_1 \le \cdots \le n_k < \infty} 1_{\{n_j = n\}} \hat{p}_{k,j}(n_1, \ldots, n_k),$$

yield the distribution of ages at which immortal mothers give birth to immortal daughters. The corresponding mean

$$\beta = \sum_{n=1}^{\infty} n e^{-\alpha n} m_n \tag{8}$$

$$= \sum_{n=1}^{\infty} E_0\left(\sum_{j=1}^{\nu} \tau_j e^{-\alpha \tau_j} 1_{\{\tau_j = n\}}\right) = E_0\left(\sum_{j=1}^{\nu} \tau_j e^{-\alpha \tau_j}\right)$$

is what is called the *average age at childbearing* (see the next section).

The factor $e^{-\alpha n_j}$ appearing in (6) reflects the reproductive value (see [6]) of the corresponding daughter. In the supercritical case $m > 1$, the population size grows exponentially at rate $\alpha > 0$, if the population escapes extinction. Thus, earlier daughters have higher reproductive value since they can contribute more to the population growth within a certain amount of time. If we compare different reproduction laws with the same marginal distribution of the offspring number $\nu$, then, clearly, the GW reproduction gives the fastest growth. In terms of the *offspring reproductive value*

$$\xi = e^{-\alpha \tau_1} + \cdots + e^{-\alpha \tau_\nu},$$

the defining equation (7) becomes $E_0(\xi) = 1$. Observe that the joint distribution of the offspring number $k$ and immortal daughter label $j$ is

$$\hat{p}_{k,j} = \sum_{1 \le n_1 \le \cdots \le n_k < \infty} \hat{p}_{k,j}(n_1, \ldots, n_k)$$

$$= 1_{\{1 \le j \le k\}} \sum_{1 \le n_1 \le \cdots \le n_k < \infty} e^{-\alpha n_j} p_k(n_1, \ldots, n_k)$$

$$= E_0(e^{-\alpha \tau_j}; \nu = k)$$

$$= p_k E_0(e^{-\alpha \tau_j} | \nu = k)$$

and the marginal size-biased distribution of the offspring number can be written as

$$\hat{p}_k = \sum_{j=1}^{k} \hat{p}_{k,j} = \sum_{j=1}^{k} \mathrm{E}_0(\mathrm{e}^{-\alpha\tau_j}; \nu = k) = \mathrm{E}_0(\xi; \nu = k) = p_k \mathrm{E}_0(\xi | \nu = k).$$

The marginal distribution of the immortal daughter label $j$ is

$$\tilde{p}_j = \mathrm{E}_0(\mathrm{e}^{-\alpha\tau_j}; \nu \geq j).$$

## 4. Size-biased CMJ processes

The immortal lineage requires an enhanced individual life space $(\hat{\Omega}_0, \hat{\mathcal{A}}_0)$ with the extended countable set of outcomes $\hat{\Omega}_0 = \bigcup_{k=0}^{\infty} \{0, 1, \ldots, k\} \times \mathbb{N}_{\uparrow}^k$ and $\hat{\mathcal{A}}_0$ being the algebra of all subsets of $\hat{\Omega}_0$. An element $\hat{\omega}_0 \in \hat{\Omega}_0$ is a vector $\hat{\omega}_0 = (j, \omega_0)$, where $\omega_0 \in \Omega_0$, as before, describes reproduction, while the additional component $j$, when positive, is meant to indicate the immortal daughter of an immortal mother. The label $j = 0$ says that the corresponding individual is mortal. We introduce two probability measures, the mortal $P_0$ and the immortal $\hat{P}_0$, on $(\hat{\Omega}_0, \hat{\mathcal{A}}_0)$, by putting

$$P_0(\hat{\omega}_0) = 1_{\{j=0\}} p_k(n_1, \ldots, n_k),$$

$$\hat{P}_0(\hat{\omega}_0) = \hat{p}_{k,j}(n_1, \ldots, n_k)$$

for $\hat{\omega}_0 = (j, n_1, \ldots, n_k)$.

With $\hat{\omega}_0 = (j, \omega_0)$, we define a new random variable by $\gamma = \gamma(\hat{\omega}_0) = j$ and update the earlier definitions by $\nu = \nu(\hat{\omega}_0) = \nu(\omega_0)$ and $\tau_i = \tau_i(\hat{\omega}_0) = \tau_i(\omega_0)$. According to the previous section,

$$\hat{P}_0(\nu = k) = \mathrm{E}_0(\xi \cdot 1_{\{\nu=k\}}),$$

as well (see (8))

$$\beta = \hat{\mathrm{E}}_0(\tau_\gamma) = \mathrm{E}_0(\tau_1 \mathrm{e}^{-\alpha\tau_1} + \cdots + \tau_\nu \mathrm{e}^{-\alpha\tau_\nu}),$$

where $\tau_\gamma$ is the regeneration age of the immortal individual with the distribution

$$\hat{P}_0(\tau_\gamma = n) = \mathrm{e}^{-\alpha n} m_n.$$

The random variable $\tau_\gamma$ is the age at childbearing, this justifying the "name mean age at childbearing" for its expectation $\beta$. The strict meaning of the age at childbearing is asymptotic ([10] and [6]). Indeed, if we sample an individual in an old supercritical population and follow her lineage backward in time, then the birth times of the individuals in the lineage form (asymptotically as the population age grows to infinity) a renewal process with the interarrival times distributed as $\tau_\gamma$.

We can now apply the approach from [5] to construct multitype branching processes based on ordering the space of individuals $\mathbb{I}$ in such a way that a mother always precedes her daughters. Note that our case corresponds to two types of individuals – mortals (type zero) and immortals (type one). Consider the product

$$(\hat{\Omega}, \hat{\mathcal{A}}) = \prod_{x \in \mathbb{I}} (\hat{\Omega}_x, \hat{\mathcal{A}}_x)$$

of enhanced life spaces $(\hat{\Omega}_x, \hat{\mathcal{A}}_x)$ which are copies of $(\hat{\Omega}_0, \hat{\mathcal{A}}_0)$. For any $x \in \mathbb{I}$, we may define the type of $x = (x_1, \ldots, x_j)$ recursively by

$$\rho_0(\hat{\omega}_0) = 1, \qquad \rho_x(\{\hat{\omega}_{x_1 \ldots x_i}\}_{i=0}^j) = 1_{\{\gamma(\hat{\omega}_{x_1 \ldots x_{j-1}}) = x_j\}}.$$

The joint distribution $\hat{P}$ for the whole size-biased CMJ process is then given by Ionescu Tulcea's theorem [11] and the transition probability for the life law of $x$, given the lives of individuals preceding $x$, being

$$P_0(\hat{\omega}_x) 1_{\{\rho_x = 0\}} + \hat{P}_0(\hat{\omega}_x) 1_{\{\rho_x = 1\}}.$$

The idea of this construction is to distinguish an immortal lineage $Y = \{y_n\}_{n \geq 0}$ with $\rho_{y_n} = 1$, where $y_0 = 0$, $y_1 = (j_0)$, $y_n = (j_0, j_{y_1}, \ldots, j_{y_{n-1}})$, given the ancestor's life $\hat{\omega}_0 = (j_0, \omega_0)$ and the lives of immortal descendents $\hat{\omega}_{y_n} = (j_{y_n}, \omega_{y_n})$. An element $\hat{\omega}_x = (j_x, \omega_x)$ of $\hat{\Omega}_x$ describes a potential branch of a tree with an immortal lineage. An element $\hat{\omega} = \{\hat{\omega}_x\}_{x \in \mathbb{I}}$ of $\hat{\Omega}$ contains the reproduction information $\omega = \{\omega_x\}_{x \in \mathbb{I}}$, along with the label set $\{j_x\}_{x \in \mathbb{I}}$. With the modified definition of the birth times $\sigma_x = \sigma_x(\hat{\omega}) = \sigma_x(\omega)$ and the tree $T = T(\hat{\omega}) = T(\omega)$, let $\hat{T} = \hat{T}(\hat{\omega}) = (T, Y)$ denote a tree with the immortal lineage. The measure $\hat{P}$ generates a size-biased measure on tree space $\hat{P}(T = t)$ as a marginal distribution of the joint distribution of $(T, Y)$.

**Example 1.** The GW case corresponds to $P(\tau_\nu \leq 1) = 1$. Then, (6) becomes

$$\hat{p}_{k,j}(1, \ldots, 1) = 1_{\{1 \leq j \leq k\}} p_k m^{-1},$$

implying (4) since $m_n = m 1_{\{n=1\}}$ and thus (7) turns to $e^{-\alpha} m = 1$. The offspring reproductive value is then simply $\xi = \nu/m$.

The beauty of the GW model lies in it being defined by a single probability distribution – that of the offspring number $\nu$. Next, we suggest two more population models (which we call *longitudinal* and *delayed* GW processes) based on a single distribution. Both assume that individual lifespan equals the offspring number $\nu$, but birth times are distributed differently. Their Malthusian parameters $\alpha_2$ and $\alpha_3$ should be compared with the GW Malthusian parameter $\alpha_1 = \ln m$.

**Example 2.** We call a CMJ process a *longitudinal GW* process if an individual lives $\nu$ years and produces exactly one offspring per year. Notationally,

$$p_k = p_k(1, \ldots, k), \qquad k \geq 1.$$

Then,

$$m_n = P(\nu \geq n)$$

and

$$\xi = e^{-\alpha} + e^{-2\alpha} + \cdots + e^{-\nu\alpha} = \frac{e^{-\alpha} - e^{-\alpha(\nu+1)}}{1 - e^{-\alpha}}$$

where $x = e^{-\alpha_2}$ solves the equation

$$E(x^\nu) = 2 - \frac{1}{x}.$$

The convex function $E(x^\nu)$ and concave function $2 - \frac{1}{x}$ intersect twice, unless $m = 1$. In the non-critical case, the relevant solution $x_2 = e^{-\alpha_2}$ is the one different from $x_1 = 1$. Clearly, $\alpha_2 < \alpha_1$. The size-biased distribution (6) becomes

$$\hat{p}_{k,j}(1, \ldots, k) = 1_{\{1 \leq j \leq k\}} e^{-\alpha_2 j} p_k$$

so that

$$\hat{p}_k = \frac{x_2 - x_2^{k+1}}{1 - x_2} p_k.$$

The longitudinal process was introduced as a mathematical artefact, for reasons of simplicity. The reader should note, though, that it generalizes to the *litter* process, where individuals live a random number of years, giving birth each year to i.i.d. litters – or, if you prefer realism to beauty, the litters follow different distributions for the first and, possibly, last years of life. Litter processes are commonly encountered in the animal world.

**Example 3.** Sevastyanov branching processes [14] are particular cases of CMJ processes where all offspring appear at the moment of their mother's death: $\tau_1 = \cdots = \tau_\nu = \tau$. The lifespan $\tau$ and the offspring number $\nu$ may depend on each other. In this case, $m_n = E(\nu; \tau = n)$ and

$$\xi = \nu e^{-\alpha\tau}, \tag{9}$$

where $\alpha$ is defined by $E(\nu e^{-\alpha\tau}) = 1$. (In the Bellman–Harris case, where $\tau$ and $\nu$ are independent, $E(e^{-\alpha\tau}) = 1/m$.)

The *delayed GW process* is a Sevastyanov process with $\tau = \nu$. As in the longitudinal GW process, an individual lives $\nu$ units of time, but, now, all of the offspring are born at the end of life. In this case,

$$E(\nu e^{-\alpha_3 \nu}) = 1.$$

We have $\alpha_3 < \alpha_2 < \alpha_1$. The size-biased distribution (6) becomes

$$\hat{p}_{k,j}(k, \ldots, k) = 1_{\{1 \leq j \leq k\}} e^{-\alpha_3 k} p_k,$$

implying that

$$\hat{p}_k = \mathrm{e}^{-\alpha_3 k} k p_k.$$

## 5. The Radon–Nikodym derivative

Let $n$ be some fixed tree level. In this section, we introduce a stopped random tree $[T]_n$ and compare its distributions under the measures $P$ and $\hat{P}$ defined earlier on the tree space. The relevant stopping procedure depends on a certain order in the CMJ tree growing algorithm. First, we determine the number of vertices and their positions on the initial branch. We then let the vertices at level $n = 1$ develop their branches. The vertices at level $n = 2$ then produce branches, and so on. The order in which the vertices at the same level develop branches is not important and is assumed to be arbitrary. The stopped tree $[T]_n$ is then obtained by stopping growth after all of the vertices at level $n$ have produced their branches. The stopped tree corresponding to an element $\omega = \{\omega_x\}_{x \in \mathbb{I}} \in \Omega$ can then be defined as

$$[T]_n = \{\omega_x, x \colon \sigma_x \le n\}.$$

For an illustration, see the tree in Figure 1, for the time being ignoring the fact that some vertices are open and others filled.

The assumed i.i.d. property of individual lives in ordinary CMJ processes leads to the following recursion:

$$P([T]_n = t) = p_k(n_1, \ldots, n_k) \prod_{i \,:\, n_i \le n} P([T]_{n-n_i} = t^{(i)}), \qquad (10)$$

where the vector $(n_1, \ldots, n_k)$ describes the root branch of the tree $t$ and $t^{(i)}$ is the daughter tree rooted at the $i$th vertex of the root branch of $t$. In Figure 1, the parameters in the recursion (10) are $k = 3$, $n_1 = 1$, $n_2 = 2$, $n_3 = 3$. There are three daughter trees in the figure: $t^{(1)}$ is to the left, $t^{(2)}$ is to the right and $t^{(3)}$ is represented by a single vertex.

**Definition 5.1.** *Let $I_n \subset \mathbb{I}$ be the set of vertices in the stopped tree $[T]_n$ which are located above the level $n$. This is called the coming generation at time $n$.*

In a similar way, we introduce a stopped tree with the immortal lineage $[\hat{T}]_n = ([T]_n, y_{u_n})$. Here,

$$u_n = \min(u \colon \sigma_{y_u} > n)$$

is the generation of the immortal lineage at time $n$ and $y_{u_n} \in I_n$ is the immortal member of the coming generation, which obviously defines the stopped immortal lineage. The stopped tree with immortal lineage has a discrete distribution satisfying a recursion similar to (10):

$$\hat{P}([\hat{T}]_n = (t, y)) = \mathrm{e}^{-\alpha n_j} p_k(n_1, \ldots, n_k) \hat{P}([\hat{T}]_{n-n_j} = (t^{(j)}, y_{-1})) \qquad (11)$$

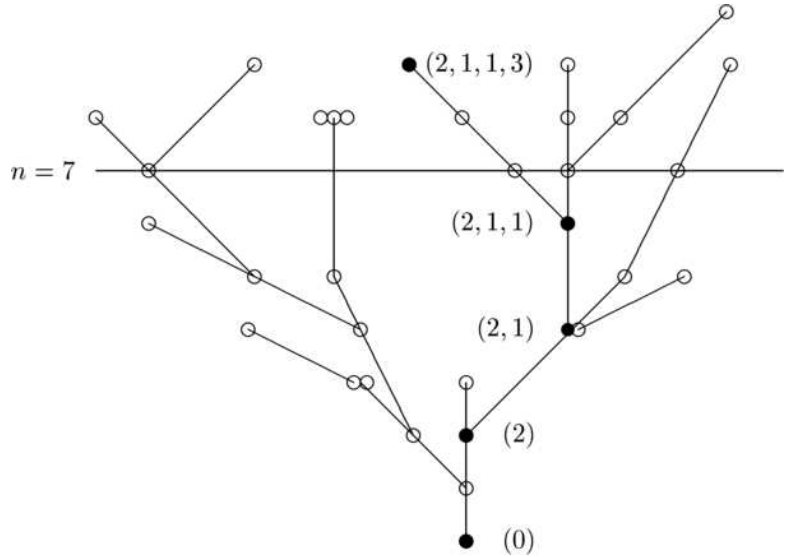

**Figure 1.** A realization of the size-biased tree $[\hat{T}]_n = (t, y)$ stopped at the level $n = 7$. The filled circles represent the immortal lineage. The coming generation at time $n$ consists of the twelve vertices above the level $n$, including the immortal vertex $y = (2, 1, 1, 3)$. In this example, the immortal lineage is stopped at generation $u_n = 4$.

$$\times \prod_{i \neq j \,:\, n_i \leq n} P([T]_{n-n_i} = t^{(i)}),$$

where $y_{-1} = (i_1, \ldots, i_u)$, given $y = (j, i_1, \ldots, i_u)$. In Figure 1, the additional parameters in the recursion (11) are $y = (2, 1, 1, 3)$, $j = 2$ and $y_{-1} = (1, 1, 3)$.

**Lemma 5.1.** *The two distributions on the tree space are related through*

$$\hat{P}([T]_n = t) = \mathrm{E}(N_n; [T]_n = t),$$

*where the Nerman martingale [9]*

$$N_n = \sum_{y \in I_n} \mathrm{e}^{-\alpha \sigma_y}$$

*is the reproductive value of the coming generation.*

**Proof.** The two recursions for the stopped trees (10) and (11) imply a recursion for the ratio

$$\frac{\hat{P}([\hat{T}]_n = (t, y))}{P([T]_n = t)} = \mathrm{e}^{-\alpha n_j} \cdot \frac{\hat{P}([\hat{T}]_{n-n_j} = (t^{(j)}, y_{-1}))}{P([T]_{n-n_j} = t^{(j)})},$$

which resolves into

$$\frac{\hat{P}([\hat{T}]_n = (t, y))}{P([T]_n = t)} = \mathrm{E}(\mathrm{e}^{-\alpha \sigma_y} | [T]_n = t).$$

It remains to recognize that, with a given tree $t$ (and therefore deterministic $I_n$),

$$\hat{P}([T]_n = t) = \sum_{y \in I_n} \hat{P}([\hat{T}]_n = (t, y))$$

$$= \sum_{y \in I_n} \mathrm{E}(\mathrm{e}^{-\alpha \sigma_y}; [T]_n = t) = \mathrm{E}(N_n; [T]_n = t). \qquad \square$$

Without loss of generality, the sequence $\{m_n\}_{n \geq 1}$ is taken as non-periodic, in that the largest common divisor of $n$ with $m_n > 0$ equals one (otherwise, we can enlarge the unit of time to ensure non-periodicity). With this restriction, Lemma 5.1 entails that $N_n = \frac{\mathrm{d}\hat{P}}{\mathrm{d}P}|_{\mathcal{B}_n}$ is the Radon–Nikodym derivative on the $\sigma$-algebra $\mathcal{B}_n$ generated by the complete lives of all individuals born up to time $n$ (even if $x$ is born at $n$, we still include her future in $\mathcal{B}_n$). It follows that the reproductive value of the coming generation $N_n$ is indeed a martingale with respect to the filtration $\mathcal{B}_n$; see [9].

Turning to our examples, we observe directly that in the GW case, $\sigma_y = n$ for all $y \in I_n$, thus confirming the corresponding result in [8] with $N_n = m^{-n} Z_n$, where $Z_n$ is the size of generation $n$. In the other two examples, despite their simple descriptions, the martingale $N_n$ does not have such an accessible form. A motivated reader would find two different expressions for $N_n$ in terms of sums of individual characteristics taken over all individuals alive at time $n$.

# 6. The $x \log x$ condition

By the martingale convergence theorem, almost surely,

$$N_n \to W \qquad \text{as } n \to \infty.$$

In the supercritical case (cf. [6]) the limit is not degenerate $P(W > 0) > 0$ if and only if

$$\mathrm{E}_0(\xi \log \xi) < \infty. \tag{12}$$

Furthermore, if we consider a function $\chi(n, w_0)$ defined on $(\Omega_0, \mathcal{A}_0)$ and put $\chi_x(n) = \chi(n, w_x)$, then the population sum

$$X_n = \sum_{x \in \mathbb{I}, \sigma_x \leq n} \chi_x(n - \sigma_x) \tag{13}$$

satisfies

$$\mathrm{e}^{-\alpha n} X_n \to \bar{\chi} \beta^{-1} W, \qquad n \to \infty, \tag{14}$$

where $\bar{\chi} = \sum_{k=0}^{\infty} e^{-\alpha k} E_0(\chi(k))$. In particular, if $\chi(n) = 1_{\{n \geq 0\}}$, then $X_n$ gives the total number of individuals born up to time $n$ and (14) holds with $\bar{\chi} = 1 - e^{-\alpha}$. Thus, in the Malthusian case (where there is a Malthusian parameter), given that the mean age at childbearing $\beta$ is also finite, the various population sums are asymptotically proportional to $\bar{\chi}$.

Next, we prove the following result concerning the $\xi \log \xi$ condition.

**Theorem 6.1.** *For supercritical CMJ processes with $\beta < \infty$, condition (12) is equivalent to*

$$E_0(\xi \log \nu) < \infty, \tag{15}$$

*which, in turn, is strictly weaker than*

$$E_0(\nu \log \nu) < \infty. \tag{16}$$

**Proof.** Since $\nu \geq \xi$, condition (12) follows from (15), which again follows from (16). The reverse is also true, provided $P_0(\tau_\nu \leq c) = 1$ for some finite $c$ because, in this case, $\xi \geq \nu e^{-\alpha c}$.

Turning to the Sevastyanov model with (9), where $\tau_1 = \cdots = \tau_\nu = \tau$ and $\beta = E_0(\tau e^{-\alpha \tau} \nu)$, we have

$$E_0(\xi \log \xi) = E_0(e^{-\alpha \tau} \nu \log \nu) - \alpha E_0(\tau e^{-\alpha \tau} \nu).$$

Thus, in the Sevastyanov case, if $\beta < \infty$, condition (12) becomes equivalent to (15). In particular, for the delayed GW process with $\tau = \nu$, conditions (15) and $\beta < \infty$ are valid, even if $E_0(\nu \log \nu) = \infty$, provided $E_0(e^{-\alpha \nu} \nu^2) < \infty$.

It remains to verify that if $\beta < \infty$ and $\hat{E}_0(\log \nu) = \infty$, which is equivalent to the negation of (15), then, necessarily, $E_0(\xi \log \xi) = \infty$. This is easy–turn to the total number of individuals $X_n$ born up to time $n$, which is defined by (13) with $\chi(n) = 1_{\{n \geq 0\}}$. Given $\beta < \infty$, the size-biased process can be bounded from below by a GW process with i.i.d. immigration distributed as $\nu - 1$ under $\hat{P}_0$. Under the assumption that $\hat{E}_0(\log \nu) = \infty$, it is straightforward to modify the argument in Section 3 of [8] and to show that $e^{-cn} X_n \to 0$ for the ordinary CMJ process, whatever the value of $c > 0$. According to [6], this implies that condition (12) does not hold. $\qquad\square$

As a final remark, we point out that condition (15) also plays a key role in the subcritical case; see [3]. It would be a useful exercise to re-prove the limit theorems for subcritical CMJ processes with a modified Lyons–Pemantle–Peres approach. This would shed new light on the limit behavior of the subcritical CMJ processes conditioned on non-extinction.

# Acknowledgements

We thank a referee for careful reading and helpful comments. This work has been supported by The Swedish Research Council and The Bank of Sweden Tercentenary Foundation.